\documentclass[12pt] {article}

\usepackage{amsmath,euscript}
\usepackage{amssymb}
\usepackage{amsthm}
\usepackage{enumerate}
\usepackage{amsfonts}
\usepackage{comment}

\usepackage{a4wide}
\usepackage{amssymb,amsmath,array}
\usepackage{amssymb}
\usepackage{amsthm}
\usepackage{enumerate}
\usepackage{amsfonts}
\usepackage{comment}

\usepackage{a4wide}
\usepackage{amssymb,amsmath,array}

\newtheorem{thm}{Theorem}[section]
\newtheorem{lem}[thm]{Lemma}

\newtheorem{prop}[thm]{Proposition}

\newtheorem{conj}[thm]{Conjecture}
\theoremstyle{defn}
\newtheorem*{defn}{Definition}

\newtheorem{rem}{Remark}
\newtheorem*{ex}{Example}

\usepackage[dvips]{graphicx}
\usepackage{amssymb}
\usepackage{amsmath}
\usepackage[latin1]{inputenc}

\usepackage[pdftex]{rotating}
\usepackage[usenames]{color}

\usepackage{graphicx,color,rotating}

\long\def\symbolfootnote[#1]#2{\begingroup%
\def\thefootnote{\fnsymbol{footnote}}\footnote[#1]{#2}\endgroup}

\begin{document}

\bibliographystyle{plain}
\begin{center}
{\bf \large Primitivity of finitely presented monomial algebras
}\\
\vspace{7 mm}
\centerline{Jason P. Bell\symbolfootnote[1]{The first author thanks NSERC
for its generous support.}} \centerline{Department of Mathematics}
\centerline{ Simon Fraser University} \centerline{ 8888 University
Dr.} \centerline{Burnaby, BC V5A 1S6.} \centerline{ CANADA}
\centerline{\tt jpb@math.sfu.ca} \vskip 5mm \centerline{Pinar
Pekcagliyan} \centerline{Department of Mathematics} \centerline{
Simon Fraser University} \centerline{ 8888 University Dr.}
\centerline{Burnaby, BC V5A 1S6.} \centerline{ CANADA}
\centerline{\tt ppekcagl@sfu.ca}\vskip 5mm \centerline{AMS Subject
Classification: 16P90} \centerline{Keywords: monomial algebras,
automata, primitivity, GK dimension, polynomial identities.}
\end{center}

\begin{abstract}
We study prime monomial algebras.  Our main result is that a prime
finitely presented monomial algebra is either primitive or it has GK
dimension one and satisfies a polynomial identity.  More generally,
we show this result holds for the class of \emph{automaton
algebras}; that is, monomial algebras that have a basis consisting
of the set of words recognized by some finite state automaton.  This
proves a special case of a conjecture of the first author and Agata
Smoktunowicz.\end{abstract}

\section{Introduction}
We consider prime monomial algebras over a field $k$.  Given a
field $k$, a finitely generated $k$-algebra $A$ is called a
\emph{monomial algebra} if $$A\ \cong \ k\{x_1,\ldots ,x_d\}/I$$
for some ideal $I$ generated by monomials in $x_1,\ldots ,x_d$.
Monomial algebras are useful for many reasons.  First, Gr\"obner
bases associate a monomial algebra to a finitely generated algebra,
and for this reason monomial algebras can be used to answer
questions about ideal membership and the Hilbert series for general
algebras.  Also, many difficult questions for algebras reduce to
combinatorial problems for monomial algebras and can be studied in
terms of forbidden subwords.  Monomial algebras are consequently a
rich area of study.  The paper of Belov,  Borisenko, and Latyshev
\cite{Belov}  is an interesting survey of what is known about
monomial algebras.

The first author and Smoktunowicz \cite{BS} studied prime monomial
algebras of quadratic growth, showing that they are either primitive
or have nonzero Jacobson radical.  By the Jacobson density theorem,
primitive algebras are dense subrings of endomorphism rings of a
module over a division algebra. For this reason, primitive ideals
are an important object of study and their study is often an
important intermediate step in classifying finite dimensional
representations of an algebra.  The first author and Smoktunowicz
\cite{BS} also made the following more general conjecture about
monomial algebras.
\begin{conj} \label{conj: 1} Let $A$ be a prime monomial algebra over a field $k$.  Then $A$ is either PI, primitive, or has nonzero Jacobson radical.
\end{conj}
These types of ``trichotomies'' are abundant in ring theory and
there are many examples of algebras for which either such a
trichotomy is known to hold or is conjectured; e.g., just infinite
algebras over an uncountable field \cite{FS},  Small's conjecture
\cite[Question 3.2]{Bell}: a finitely generated prime Noetherian
algebra of quadratic growth is either primitive or PI.

We prove Conjecture \ref{conj: 1} in the case that $A$ is finitely
presented; in fact we are able to prove it more generally when $A$
is a prime \emph{automaton algebra}; that is, when $A$ is a monomial
algebra with a basis given precisely by the words recognized by a
finite state machine.

Our main result is the following theorem.
\begin{thm} \label{thm: main} Let $k$ be a field and let $A$ be a
prime finitely presented monomial $k$-algebra.  Then $A$ is either
primitive or $A$ satisfies a polynomial identity.
\end{thm}

A consequence of this theorem is that any finitely generated prime monomial ideal $P$ in a finitely generated free algebra $A$ is necessarily primitive, unless $A/P$ has GK dimension at most $1$.

We prove our main result by showing that a finitely presented monomial algebra $A$ has a well-behaved free subalgebra.  In this case, well-behaved means that there the poset of left ideals of the subalgebra embeds in the poset of left ideals of $A$ and nonzero ideals in our algebra $A$ intersect the subalgebra non-trivially.  A free algebra is primitive if it is free on at least two generators; otherwise it is PI.  From this fact and the fact that $A$ has a well-behaved free subalgebra, we are able to deduce that $A$ is either primitive or PI.

In Section $2$ we give some background on finite state automata and automaton algebras.  In Section $3$ we give some useful facts about primitive algebras and PI algebras that we use in obtaining our dichotomy.  In Section $4$ we prove our main result.
\section{Automata theory and automaton algebras}

In this section we give some basic background about finite state automata and automaton algebras.
A finite state automaton is a machine that accepts as input words on a finite alphabet $\Sigma$ and has a finite number of possible outputs.  We give a more formal definition.

\begin{defn}{\em
A \emph{finite state automaton} $\Gamma$ is a 5-tuple $(Q, \Sigma , \delta, q_0, F)$,
where:
\begin{enumerate}
\item $Q$ is a finite set of states;
\item $\Sigma$ is a finite alphabet;
\item $\delta : Q \times \Sigma \rightarrow Q$ is a transition
function;
\item $q_0 \in Q$ is the initial state;
\item $F \subseteq Q$ is the set of accepting states.
\end{enumerate} }
\end{defn}
We refer the reader to Sipser \cite[page 35]{Sipser} for more
background on automata. We note that we can inductively extend the
transition function $\delta$ to a function from $Q\times \Sigma^*$
to $Q$, where $\Sigma^*$ denotes the collection of finite words of
$\Sigma$.\footnote{We simply define $\delta(q, \varepsilon)=q$ if
$\varepsilon$ is the empty word and if we have defined $\delta(q,w)$
and $x\in \Sigma$, we define $\delta(q,wx):=\delta(\delta(q,w),x)$.}

We give an example.
\begin{ex} {\em The finite state machine described in Figure 1 has
alphabet $\Sigma = \{0,1\}$, states $\{q_0,q_1\}$ accepting state
$\{q_0\}$ and transition rules $\delta(q_i,1)=q_{1-i}$,
$\delta(q_i,0)=q_i$ for $i=0,1$.  In particular, if $w$ is a word on
$\{0,1\}$ then $\delta(q_0,w)$ is $q_0$ if and only if the number of
ones in $w$ is even.}
\end{ex}

\begin{figure}
\begin{center}
\includegraphics{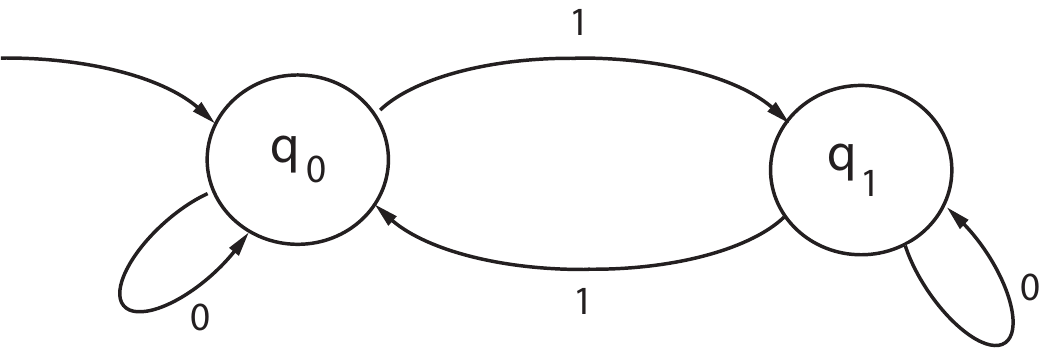}
\caption{A finite state automaton with two states.}
\end{center}
\end{figure}

\begin{defn} {\em Let  $\Gamma=(Q, \Sigma , \delta, q_0, F)$
be a finite state automaton.  We say that a word $w\in \Sigma^*$ is
\emph{accepted} by $\Gamma$ if $\delta(q_0,w)\in F$; otherwise, we
say $w$ is \emph{rejected} by $\Gamma$.}
\end{defn}

We now describe the connection between monomial algebras and finite
state automata.

\begin{defn}{\em Let $k$ be a field and let $A=k \{ x_1, x_2, \dots, x_n \} /
I$ be a monomial algebra. We say that $A$ is an \emph{automaton
algebra} if there exists a finite state automaton $\Gamma$ with
alphabet $\Sigma = \{x_1,\ldots ,x_n\}$ such that the words $w$ is
accepted by $\Gamma$ if and only if $w\not\in I$. }
\end{defn}
\begin{rem}
Any finitely presented monomial algebra is an automaton algebra.
\end{rem}
\noindent {\bf Proof.} See Belov, Borisenko, and
Latyshev \cite[Proposition 5.4 p. 3528]{Belov}. \qed
\vskip 2mm

We note that in order for finite state automaton $\Gamma$ to give
rise to a monomial algebra, the collection of words that are
rejected by $\Gamma$ must generate a two-sided ideal that does not
contain any two sided words. This need not occur in general (see,
for example, Figure 1, in which we take $F=\{q_0\}$).

In general, an automaton algebra can have many different corresponding finite state automata.  We may assume, however, that the corresponding finite state automaton is \emph{minimal}.

\begin{defn} {\em We say that a finite state automaton  $\Gamma=(Q, \Sigma , \delta, q_0, F)$ is \emph{minimal} if for $q_1,q_2\in Q$ with $q_1\not =q_2$ we have
\[ \left\{ w\in \Sigma^*~:~\delta(q_1,w)\in F\right\} \not =  \left\{ w\in \Sigma^*~:~\delta(q_2,w)\in F\right\}\]
and for every $q\in Q$ there is a word $w\in \Sigma^*$ such that $\delta(q_0,w)=q$.
}
\end{defn}
We note that this definition of minimality is slightly different from other definitions that appear in the literature.  It can, however, be shown to be equivalent \cite{W1}.
By the Myhill-Nerode theorem \cite{HU}, if $A$ is an automaton algebra, there is a minimal automaton $\Gamma$ such that the images in $A$ of the words accepted by $\Gamma$ form a basis for $A$; moreover, this automaton $\Gamma$ is unique up to isomorphism.  For this reason, we will often speak of \emph{the} minimal automaton corresponding to an automaton algebra. 

\section{Primitivity and Polynomial identities}
In this section we give some important background on primitive
algebras and algebras satisfying a polynomial identity.  We first
recall the definitions of the two main concepts that make up our
dichotomy.
\begin{defn}
{\em A ring $R$ is \emph{left-primitive} if it has a faithful simple
left $R$-module $M$. }
\end{defn}

Right-primitivity is defined analogously.  Left-primitivity and
right-primitivity often coincide; nevertheless there are examples of
algebras which are left- but not right-primitive \cite{Bergman}. For
the purposes of this paper, we will say that an algebra is
\emph{primitive} if it is \emph{both} left- and right-primitive.

Two useful criteria for being primitive are given in the following remark.

\begin{rem} Let $R$ be a ring with unit.  The following are equivalent:
\begin{enumerate}
\item $R$ is left primitive;
\item $R$ has a maximal left
ideal $I$ that does not contain a nonzero two sided ideal of $R$;
\item $R$ has a left ideal $I$ such that $I+P=R$ for every nonzero prime ideal $P$.
\end{enumerate}
\end{rem}
\noindent {\bf Proof.} For the equivalence of (1) and (2), see Rowen
\cite[Page 152]{Rowen}; for the equivalence of (2) and (3) note that
(2) trivially implies (3), and if $I$ is a left ideal that does not
contain a nonzero two-sided ideal of $R$ then by Zorn's lemma we can
find a maximal left ideal with this property. \qed \vskip 2mm The
other concept used in the dichotomy in the statement of Theorem
\ref{thm: main} is that of being \emph{PI}.
\begin{defn}
{\em We say that a $k$-algebra $A$ satisfies a \emph{polynomial
identity} if there is a nonzero noncommutative polynomial
$p(x_1,\dots ,x_n) \in k\{x_1,\dots,x_n\}$ such that
$p(a_1,\dots,a_n)=0$ for all $(a_1,\dots,a_n) \in A^n$.  If an
algebra $A$ satisfies a polynomial identity we will say that $A$ is
\emph{PI}.}
\end{defn}
Polynomial identity algebras are a natural generalization of
commutative algebras, which, by definition, satisfy the polynomial
identity $xy-yx=0$.  An important theorem of Kaplansky \cite[p. 157, 6.3.1]{Herstein} shows that an algebra that is \emph{both}
primitive and PI is a matrix ring over a division algebra that is
finite dimensional over its centre. Kaplansky's theorem shows that
being primitive and being PI are in some sense incongruous and this
incongruity is expressed in the fact that for many classes of
algebras there are either theorems or conjectured dichotomies which
state that the algebra must be either primitive or PI \cite{FS, BS,
Bell}.

To prove Theorem \ref{thm: main}, we rely on a reduction to free
algebras.  A free algebra on $1$ generator is just a polynomial ring
and hence PI.  A free algebra on two or more generators is
necessarily primitive.  This result is due to Samuel  \cite[page
36]{Jacobson}. We outline a proof of this fact, since Jacobson only
describes the result for free algebras on precisely two generators.

\begin{thm} \label{Jacobs}
A free algebra that is either countably infinitely generated or is
generated by $d$ elements for some natural number $d\ge 2$ is
primitive.
\end{thm}
\noindent {\bf Proof.}  Since a free algebra is isomorphic to its
opposite ring, it is sufficient to prove left primitivity. First, if
$A$ is the free algebra on two generators, say $A=k\{x,y\}$, then we
construct a left $A$-module $M$ as follows.  Let $M=\sum_{i\ge 0}
ke_i$ and let $A$ act on $M$ via the rules
\[ xe_i \ = \ e_{i-1} \qquad {\rm and} \qquad ye_i = e_{i^2+1}, \]
where we take $e_{-1}=0$.   Then $M$ is a faithful simple left
$A$-module and so $A$ is left primitive.  Next observe that if
$A=k\{x,y\}$ then $k+Ay$ is free on infinitely many generators
$y,xy,x^2y,\ldots $ and hence a free algebra on a countably infinite
number of generators is primitive \cite{LRS}.  It follows that if
$d\ge 2$ is a natural number then $A=k\{x_1,\ldots ,x_d\}$ is again
primitive since $k+Ax_d$ is free on a countably infinite number of
generators \cite{LRS}.   \qed \vskip 2mm We use Theorem \ref{Jacobs}
to prove the primitivity of prime automaton algebras of GK dimension
greater than $1$.  To do this, we use a result that allows one to
show that an algebra with a sufficiently well-behaved primitive
subalgebra is itself necessarily primitive; we make this more
precise in Proposition \ref{prop: nearly free}, but before stating
this proposition we require a definition.
\begin{defn}
{\em Let $B$ be a subring of a ring $A$.  We say that $A$ is
\emph{nearly free} as a left $B$-module if there exists some set
$E=\{x_{\alpha}~|~\alpha\in S\}\subseteq A$ such that:
\begin{enumerate}
\item $x_{\alpha_0}=1$ for some $\alpha_0\in S$;
\item $A=\sum_{\alpha} Bx_{\alpha};$
\item  if $b_1,\ldots ,b_n\in B$ and $b_1x_{\alpha_1}+\dots +b_nx_{\alpha_n}=0$
then $b_ix_{\alpha_i}=0$ for $i=1,\ldots,n$.

\end{enumerate} }
\end{defn}
We note that it is possible to be nearly free over a subalgebra
without being free.  For example, let $A=\mathbb{C} [x] /(x^3)$ and
let $B$ be the subalgebra of $A$ generated by the image of $x^2$ in
$A$.  Then $A$ is $3$-dimensional as a $\mathbb{C}$-vector space
while $B$ is $2$-dimensional.  Hence $A$ cannot be free as a left
$B$-module.  Let $\overline{x}$ denote the image of $x$ in $A$. Then
$$A \ =\  B+B\overline{x}.$$  Moreover $A$ is $\mathbb{N}$-graded
and $B$ is the graded-subalgebra generated by homogeneous elements
of even degree.  Hence if $b_1 + b_2 \overline{x}=0$ then
$b_1=b_2\overline{x}=0$.  Hence $A$ is nearly free as a $B$-module.

\begin{prop}\label{prop: nearly free} Let $A$ be a prime algebra and suppose that $B$ is a right primitive subalgebra of $A$ such that:
\begin{enumerate}
\item $A$ is nearly free as a left $B$-module;
\item every nonzero two-sided ideal $I$ of $A$ has the property that $I\cap B$ is nonzero.
\end{enumerate}
Then $A$ is right primitive.
\end{prop}
\noindent {\bf Proof.}  Pick a maximal right $I$ of $B$ that does not contain a nonzero two-sided ideal of $B$.
Let $E=\{ x_\alpha : \alpha \in S \}$ be a subset of $A$ satisfying:
\begin{enumerate}
\item $A=\sum_{x_{\alpha} \in E} Bx_{\alpha}$;
\item if
$b_1x_{\alpha_1}+ \cdots +b_dx_{\alpha_d}=0$, then $b_ix_{\alpha_i}=0$ for every $i$;
\item $x_{\beta}=1$ for some $\beta\in S$.
\end{enumerate}
Then $$IA \ = \ \sum_{\alpha\in S} Ix_{\alpha}.$$ We claim that $IA$
is proper right ideal of $A$.  If not then $$1=x_\beta=\sum a_k
x_{\alpha_k},$$ for some $a_k\in I$. Since $A$ is nearly free as a
left $B$-module, $x_\beta -a x_\beta=0$ for some $a\in I$,
contradicting the fact that $I$ is proper.  Thus $IA$ is a proper
right ideal. By Zorn's lemma, we can find a maximal right ideal $L$
lying above $IA$.  We claim that $A/L$ is a faithful simple right
$A$-module.  To see this, suppose that $L$ contains a nonzero prime
ideal $P$ of $A$.  By assumption, $P\cap B=Q$ is a nonzero ideal of
$B$ and is contained in $L\cap B=I$, a contradiction.  The result
follows. \qed

\section{Proofs}
In this section we prove the following generalization of Theorem \ref{thm: main}.
\begin{thm} \label{thm: mainx} Let $k$ be a field and let $A$ be a
prime automaton algebra over $k$.  Then $A$ is either
primitive or $A$ satisfies a polynomial identity.
\end{thm}

To prove this we need a few definitions.

\begin{defn} {\em
Let $\Gamma = (Q,\Sigma, \delta, q_0, F)$ be a minimal finite state automaton. Given a state $q\in Q$, we say a word $w \in
\Sigma^*$ is $q$-\emph{revisiting} if $w=w'w''$
for some $w',w''\in \Sigma^*$ with $w'$ non-trivial such that $\delta(q,w')=q$. Otherwise, we say $w$ is
$q$-\emph{avoiding}.}
\end{defn}

A key obstruction in this proof is that there exist examples of prime automaton algebras for which  there are no non-trivial words in $\Sigma^*$ that are $q_0$-revisiting in the corresponding minimal automaton $\Gamma = (Q,\Sigma, \delta, q_0, F)$.  For example, if $A=k\{x,y\}/(x^2,y^2)$.  Then
the minimal automaton corresponding to the algebra $A$ is given in Figure
2, where the accepting states are $F=\{q_0,q_1,w_2\}$.  We note that it is impossible to revisit the initial state $q_0$ in this case.
\begin{figure}
\includegraphics{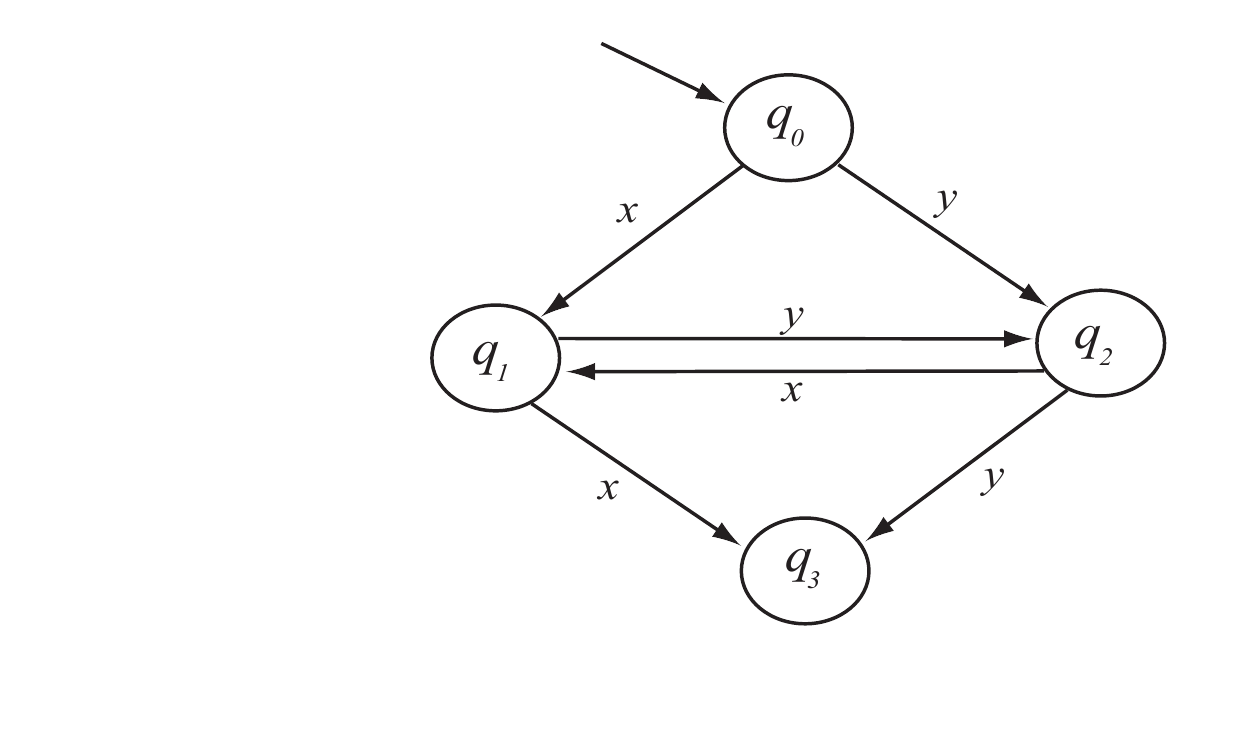}
\caption{A finite state automaton in which the only $q_0$-revisiting word is the empty word.}
\end{figure}

Before proceeding with the generalization of Theorem \ref{thm:
main}, we define an equivalence relation on the accepting states of
a minimal finite state automaton $\Gamma = (Q,\Sigma, \delta, q_0, F)$. We
say that $q_i\sim q_j$ if there exists words $w$ and $w'$ such that
$\delta(q_i,w)=q_j$ and $\delta(q_j,w')=q_i$.

Given a minimal finite state automaton $\Gamma = (Q,\Sigma, \delta, q_0,
F)$, we can put a partial order between the equivalence classes in
the following way.  Let $q$ and $q'$ be two accepting states and let
$[q]$ and $[q']$ denote their equivalence classes. We say that
$[q]\le [q']$ if there is a word $w$ such that $\delta(q,w)=q'$.
(Note that if $[q]\le [q']$ and $[q']\le [q]$ then $q\sim q'$ and so
the two classes are the same.)  Figure 3 gives an example of a finite state automaton in which the set of states has been partitioned into equivalence classes.  In this example, the class $E_1$ is minimal, $E_4$ is maximal, and $E_2$ and $E_3$ are incomparable.  

\begin{figure}
\begin{center}
\includegraphics{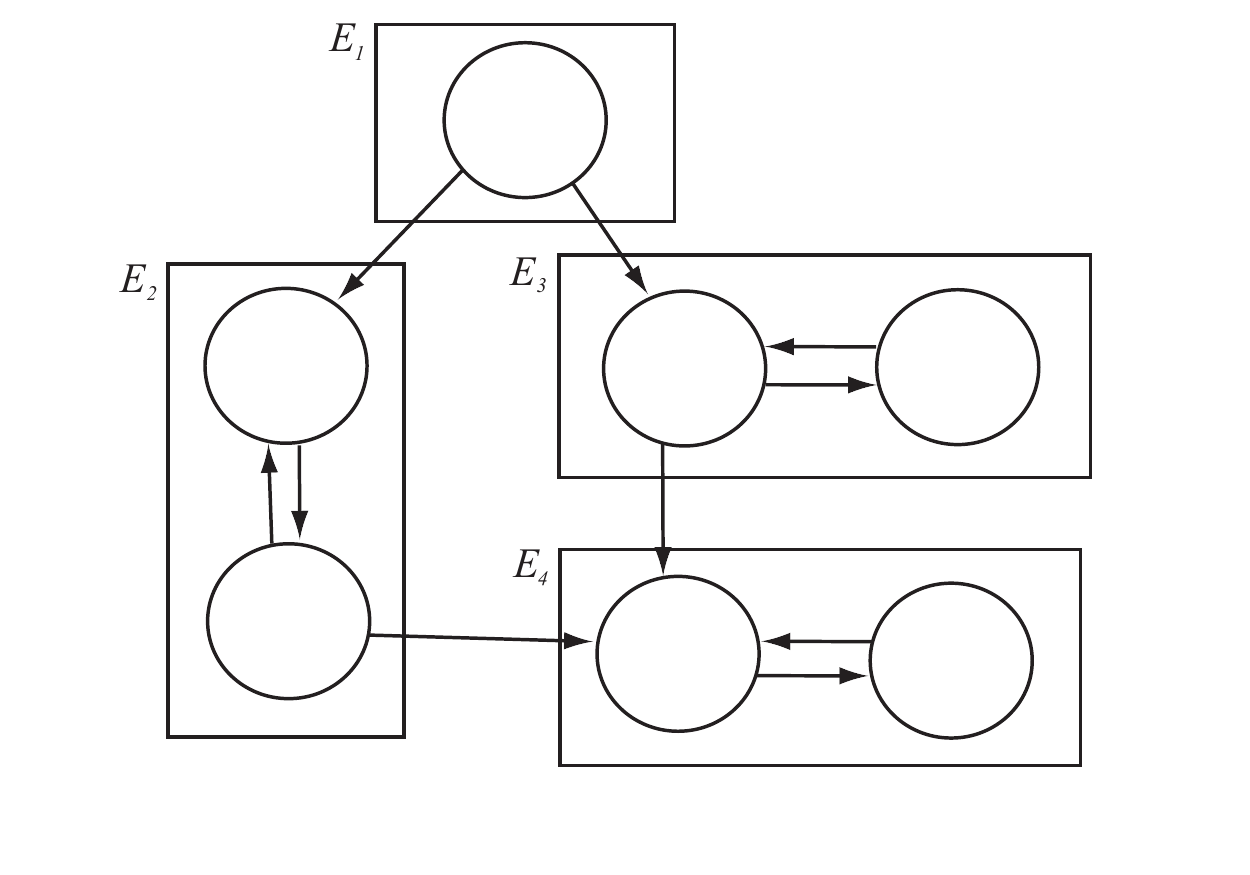}
\end{center}
\caption{A finite state machine with four equivalence classes.}
\end{figure}

To obtain the proof of Theorem \ref{thm: main}, we show that $A$ has
a well-behaved free subalgebra.  We call the subalgebras we
construct \emph{state subalgebras}.
\begin{defn} {\em
Let $A$ be an automaton algebra and let $\Gamma=(Q, \Sigma, \delta,
q_0, F)$ be its corresponding minimal finite state automaton. Given a state $q
\in F$, we define the \emph{state subalgebra of }$A$\emph{
corresponding to} $q$ to be the subalgebra generated by all words
$w\in \Sigma^*$ such that $\delta(q,w)=q$.}
\end{defn}

\begin{lem} \label{lem:state subalgebra}
Let $A$ be an automaton algebra and let $\Gamma=(Q, \Sigma, \delta, q_0,
F)$ be its corresponding minimal finite state automaton.  A state
subalgebra $B$ of $A$ corresponding to some state $q$ in $F$ is a free algebra.
\end{lem}
\noindent {\bf Proof.} We claim that $B$ is free on $$E \ := \
\left\{w\in\Sigma^* \, \left| \  \delta(q,w)=q,~{\rm every~ nonempty~ proper~
initial~ subword~ of~}w~{\rm is} ~q\rule[.025in]{.05in}{.009in}{\rm
avoiding}\right. \right\}.$$ Since $B$ is generated by words $w$
such that $\delta(q,w)=q$ and every such word $w$ can be decomposed
into a product of words $w=w_1\cdots w_d$ with $\delta(q,w_i)=q$ and
for which every nonempty proper initial subword of $w_i$ is
$q$-revisiting, we see that $B$ is generated by $E$. Suppose that
$B$ is not free on $E$. Then we have a non-trivial relation of the
form $$\sum c_{i_1,\ldots ,i_d} w_{i_1}\dots w_{i_d} \ = \ 0,$$ in
which only finitely many of the $c_{i_1,\ldots ,i_d}$ are nonzero
and each $w_{i_1},\ldots ,w_{i_d}\in E$. Since $A$ is a monomial
algebra, we infer that we must have a relation of the form
$$w_{i_1} \dots w_{i_d} \ = \  w_{j_1} \dots w_{j_e}$$ with
$$(w_{i_1},\ldots ,w_{i_d})\not = (w_{j_1},\ldots ,w_{j_e}).$$
Pick such a relation with $d$ minimal.  Then note that $w_{i_1}\not
= w_{j_1}$ for otherwise, we could remove $w_{i_1}$ from both sides
and have a smaller relation: $$w_{i_2} \dots w_{i_d} \ = \  w_{j_2}
\dots w_{j_e}.$$ But then either $w_{i_1}$ is a proper
$q$-revisiting initial subword of $w_{j_1}$ or vice versa, which is
impossible by the definition of the set $E$.  The result follows.
\qed
\begin{lem} \label{lem:ideals lift}
Let $A$ be an automaton algebra and let $\Gamma=(Q, \Sigma, \delta, q_0,
F)$ be its corresponding minimal finite state automaton.  If $B$ is a state subalgebra of $A$ corresponding to some state $q$ in $F$ then $A$ is nearly free as a left $B$-module.
\end{lem}
\noindent {\bf Proof.} Let $$E=\{1\}\cup\left\{w\in\Sigma^* \, \left| \  \delta(q_0,w)\in F~{\rm and}~w~{\rm is} ~q\rule[.025in]{.05in}{.009in}{\rm
avoiding}\right. \right\}.$$  The condition that $\delta(q_0,w)\in F$ is just saying that $w$ has nonzero image in $A$.
We claim first that
$$A \ = \ \sum_{x\in E} Bx.$$
Since $A$ is spanned by words, it is sufficient to show that every word $w$ with nonzero image in $A$ is of the form $bx$ for some $b\in B$ and $x\in E$.  Note, however, there is some proper initial subword $b$ of $w$ such that $w=bx$, $\delta(q,b)=q$ and $x$ is either $q$-avoiding of $x=1$.  Thus we obtain the first claim.

Next observe that if
$$\sum_{i=1}^d b_i x_i \ = \ 0,$$ with $x_i\in E$, $b_i\in B$ then we must have $b_1 x_1 = \cdots =b_dx_d=0$.  To see this, observe by the argument above, every word $u$ has a unique expression as $bx$ for some word $b\in B$ and $x\in E$.  Suppose that $$\sum_{i=1}^d b_i x_i \ = \ 0$$ and $b_1x_1\not =0$.  Then there is some word $u$ which appears with a nonzero coefficient in $b_1x_1$.  But by the preceding remarks, $u$ cannot appear with nonzero coefficient in any of $b_2x_2,\ldots , b_dx_d$.  Since $A$ is a monomial algebra, we obtain a contradiction.  Thus $A$ is nearly free as a left $B$-module. \qed
\vskip 2mm
We have now shown that a prime automaton algebra as a free subalgebra $B$ such that $A$ is nearly free as a left $B$-module.  To complete the proof that $A$ is primitive or PI, we must show that nonzero ideals of $A$ intersect certain state subalgebras non-trivially.

\begin{prop} \label{prop:ideal}
Let $A$ be a prime automaton algebra with corresponding minimal finite state
automaton $\Gamma=(Q,\Sigma, \delta, q_0, F)$.  Suppose $q \in F$ is in a maximal equivalence class of $F$ under the order described above and $B$ is the state subalgebra
corresponding to $q$.  If $I$ is a nonzero two sided ideal of
$A$ then $I \cap B$ is a nonzero two sided ideal of $B$.
\end{prop}
\noindent {\bf Proof.} Every element $x\in I$ can be written as
$$\sum c_w w,$$ where $w\in \Sigma^*$.  Among all nonzero $x\in I$, pick an element
$$x=c_1 w_1 + \cdots + c_d w_d$$ with $d$ minimal.  Then we may assume $c_1,\ldots ,c_d$ are all nonzero.
Pick $u$ such that $\delta(q_0,u)=q$.  Since $A$ is prime, there is some word $v$ such that $uvx\not = 0$.  Then $uvx = c_1 uvw_1 + \cdots + c_d uvw_d$ is a nonzero element of $I$.  My mimimality of $d$, $uvw_i$ has nonzero image in $A$ for every $i$.  Consequently, $\delta(q_0,uvw_i)\in F$ for all $i$.  Since $q$ is in a maximal equivalence class of $F$ and $\delta(q_0,u)=q$, $\delta(q_0,uvw_i)\in [q]$ for $1\le i\le d$.

Note that if $\delta(q_0,uvw_i)\not = \delta(q_0,uvw_j)$ for some $i,j$ then by minimality of $\Gamma$, there is some word $w\in \Sigma^*$ such that $\delta(q_0,uvw_iw)\in F$ and $\delta(q_0,uvw_jw)\not\in F$ (or vice versa).  Consequently, $uvw_iw$ has nonzero image in $A$ and $uvw_jw=0$ in $A$.  Thus
$uvxw$ is a nonzero element of $I$ with a shorter expression than that of $x$, contradicting the minimality of $d$.  It follows that $$\delta(q_0,uvw_1)=\cdots = \delta(q_0,uvw_d).$$
Since $\delta(q_0,uvw_1)\in [q]$, there is some word $u'$ such that $\delta(q_0,uvw_1u')=q$.  Consequently,
$$\delta(q_0,uvw_1u')\ = \ \cdots \ = \ \delta(q_0,uvw_du')\ = \ q.$$
Thus $uvw_iu'\in B$ for $1\le i\le d$ and so $uvxu'\in B\cap I$ is nonzero.  The result follows. \qed
\vskip 2mm

\noindent {\bf Proof of Theorem \ref{thm: mainx}.}  It is sufficient to show that $A$ is right primitive since the opposite ring of $A$ is again an automaton algebra.\footnote{By Kleene's theorem \cite[Theorem 4.1.5]{AS}, the collection of words accepted by a finite state automaton forms a regular language; by symmetry in the definition of a regular language, the \emph{reverse} language obtained by reversing all strings in a given regular language is again regular.  At the level of algebras, string reversal corresponds to multiplication of words in the opposite ring.}
Let $\Gamma=(Q, \Sigma, \delta,
q_0, F)$ be the minimal finite state automaton corresponding to $A$.
We pick a state $q\in F$ that is in an equivalence class $[q]$ that is maximal with respect to the order described above.  We let $B$ be the state subalgebra of $A$ corresponding to $q$.  By Lemma \ref{lem:state subalgebra}, $B$ is a free algebra.  We now have two cases.
\vskip 2mm
\noindent {\bf Case I:} $B$ is free on at most one generator.
\vskip 2mm
\noindent In this case, we claim that $A$ satisfies a polynomial identity.  Let $u$ be a word satisfying
$\delta(q_0,u)=q$.  Let $v$ be a word with nonzero image in $A$.  Since $A$ is prime, there is some word $w$ such that $uwv$ has nonzero image in $A$.  Thus $\delta(q_0,uwv)\in F$.  Since $\delta(q_0,u)=q$, and $[q]$ is a maximal equivalence class, $\delta(q,wv)\in [q]$.  In particular, there is some word $t$ such that $\delta(q_0,uwvt)=q$.  Thus $wvt\in B$.  But $B$ is free on at most one generator.  In particular every word in $B$ must be a power of some (possibly empty word) $b$.  Thus $v$ is a subword of $b^m$ for some $m$.  It follows that every word with nonzero image in $A$ is a subword of $b^m$ for some $m$.  In particular, the number of words in $A$  of length $n$ that have nonzero image in $A$ is bounded by the length of $b$.  Hence $A$ has GK dimension at most one (cf. Krause and Lenagan \cite[Chapter 1]{KL}). Thus $A$ is PI \cite{SW}.
\vskip 2mm
\noindent {\bf Case II:} $B$ is free on two or more generators.
\vskip 2mm
\noindent In this case, $B$ is primitive by Theorem \ref{Jacobs}.
By Lemma \ref{lem:ideals lift}, $A$ is nearly free as a left $B$-module.  By Proposition \ref{prop:ideal}, nonzero ideals of $A$ intersect $B$ non-trivially.  Hence $A$ is right primitive by Proposition \ref{prop: nearly free}.  The result follows. \qed

\end{document}